\documentclass[11pt]{article}
\usepackage[margin=1in]{geometry}  
\usepackage{graphicx}
\usepackage{amssymb}
\usepackage{amsmath}
\usepackage{caption}
\usepackage{subcaption}
\usepackage{float}
\usepackage{url}

\title{The Shape of Symmetric Binary Trees}
\author{Larry Riddle\thanks{\texttt{lriddle@agnesscott.edu}} \\
Agnes Scott College}
\date{}

\begin{document}
\maketitle
\begin{abstract}
Mandelbrot and Frame studied the geometry of self-contracting symmetric binary trees in which they stated that the height of such trees occurred at the branch tip of the path consisting of branches that alternate left and right. Taylor proved that this happens for both self-avoiding as well as self-contacting symmetric binary trees (if we ignore the height of the trunk and just consider the branch tips). In his commentary on Mandelbrot's and Frame's work, West gave an example of a self-overlapping tree in which this  alternating left-right path does not give the highest point of the tree, and said that more analysis was needed. In this paper we show how such examples can be constructed for all but a countable number of angles. We also investigate the conditions for when the sides and bottom of a self-overlapping symmetric binary tree differ from what happens with self-avoiding and self-contacting trees.
\end{abstract}
\section{Introduction}
To construct a symmetric binary tree, choose an angle $\theta$ with $0^\circ < \theta < 180^\circ$ and a scaling factor $r$ with $0 < r < 1$. Start with a vertical line segment (the trunk) of length 1. The trunk splits into two branches at the top that each form an angle of $\theta$ with the linear extension of the trunk, one to the left and one to the right. Each branch has length $r$. Each of these two branches forms the trunk of a subtree that splits into two more branches following the same rule. The angle is again $\theta$ and the length of each of the four new branches is $r^2$,
\begin{figure}[ht] 
   \centering
   \includegraphics[width=2in]{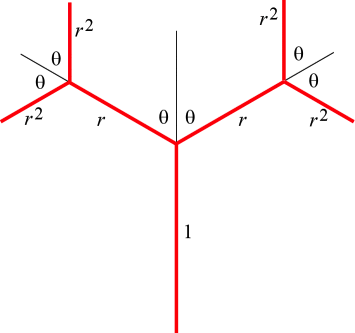} 
   \caption{Two iterations}
   \label{fig:construction}
\end{figure}
as shown in Figure~\ref{fig:construction}. The symmetric binary tree is obtained by continuing to add more branches ad infinitum, using the angle $\theta$ and scaling factor $r$ for each set of new branch segments. The limit points of the branches in a binary tree are called the branch tips. 
\begin{figure}[ht]
   \centering
   \begin{subfigure}[b]{0.3\textwidth}
       \includegraphics[width=\textwidth]{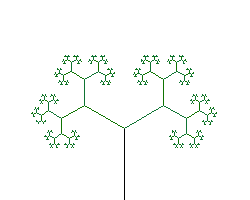}
       \caption{$\theta = 60^\circ, r = 0.58$}
   \end{subfigure}
   \begin{subfigure}[b]{0.3\textwidth}
       \includegraphics[width=\textwidth]{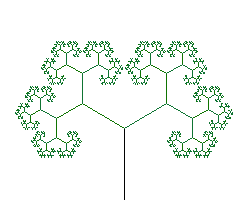}
       \caption{$\theta = 60^\circ, r = 0.61803$}
   \end{subfigure}   
   \begin{subfigure}[b]{0.3\textwidth}
       \includegraphics[width=\textwidth]{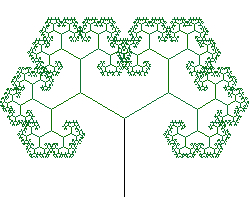}
       \caption{$\theta = 60^\circ, r = 0.65$}
   \end{subfigure}
   \caption{Symmetric binary trees}
   \label{fig:examples}
\end{figure}
Figure~\ref{fig:examples} illustrates examples of three symmetric binary trees. If the scaling factor $r$  is too small, the branches of the tree will be self-avoiding, while if $r$ is too large the branches will overlap. Mandelbrot and Frame [1,2] showed that for each angle $\theta$, there is a unique scaling factor $r_{sc}$ such that the symmetric binary tree with that $\theta$ and $r = r_{sc}$ will be self-contacting, that is, branches may touch at a single point but may not cross. In this case, the branch tips of some left-side branches may coincide with branch tips of some right-side branches, but no branch tip will coincide with any non-tip point of the tree. In Figure 2, the three trees are self-avoiding, self-contacting, and self-overlapping, respectively.

The paths of a symmetric binary tree can be addressed by finite strings of the letters $L$ and $R$, with $L$ corresponding to branching to the left and $R$ corresponding to branching to the right. The branch tips are obtained as infinite sequences of $L$ and $R$. For example, the sequence $(RL)^\infty = RLRLRL\dots$ would be a path alternating right and left branches. 
\begin{figure}[ht] 
   \centering
   \includegraphics[width=3.5in]{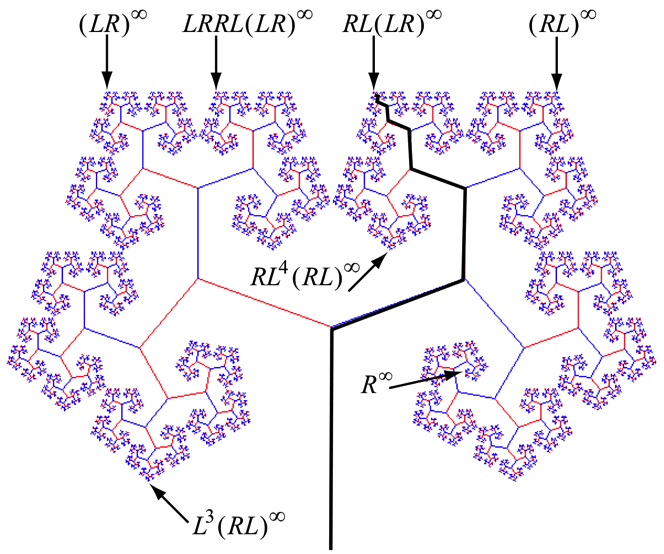} 
   \caption{Branch tip labeling}
   \label{fig:tip labeling}
\end{figure}
In Figure~\ref{fig:tip labeling}, notice the symmetry with the address $(LR)^\infty$. The address $R^\infty$ corresponds to an infinite spiral of right branches.

\begin{figure}[ht] 
   \centering
   \includegraphics[width=4in]{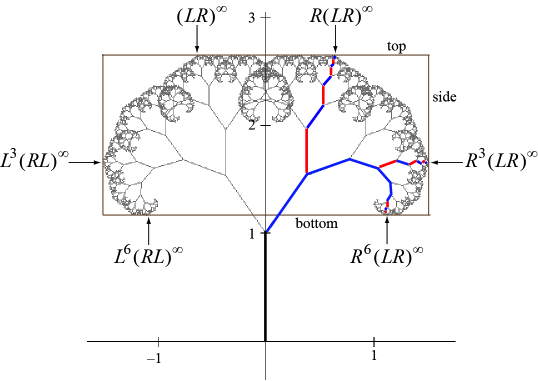} 
   \caption{$\theta = 35^\circ, r = 0.6$ (top = 2.6539, side = 1.51964, bottom = 1.15588)}
   \label{fig:dimensions}
\end{figure}

We will take the base of the tree to be at the origin and the initial vertical trunk to be of length 1. By the ``top" of the tree we will mean the top of the branch tips, that is, the largest y-coordinate of the branch tips. The ``bottom" of the tree will mean the bottom of the branch tips, that is, the smallest $y$-coordinate of the branch tips. The ``right side" of the tree will be the largest $x$-coordinate of the branch tips (and by symmetry, the left side will be the negative of this $x$-coordinate.) The trunk is not considered in determining the top or bottom of the tree. Figure~\ref{fig:dimensions} illustrates these dimensions for a tree with $\theta = 35^\circ$ and $r = 0.6$, along with some addresses of  branch tips at which the extreme $y$ and $x$-coordinates are obtained. 

Details for many of the results and computations that follow may be found at the website \url{https://larryriddle.agnesscott.org/ifs/pythagorean/symbinarytreeShape.htm}.

\section{Top of the Tree}
To get to the top of the tree, we want to go vertically as much as possible for each new branch in the construction of the tree. This suggests that we want to alternate left and right branches and follow the path $(LR)^\infty$, or by symmetry, the path $(RL)^\infty$. Mandelbrot and Frame first discussed the height of a self-contacting tree in [1]. Taylor [4] proved that for self-contacting and self-avoiding trees, the $y$-coordinate of the point with address $(LR)^\infty$ does indeed give the top of the tree. There are, in fact, many branch tips with the same $y$-coordinate as that for the $(LR)^\infty$ path, such as the branch tip for $(RL)^\infty$. Consider the finite path with address $P_n = T_1 T_2 T_3 \dots T_n$ where $T_i$ is either the pair $LR$ or the pair $RL$. Then the $y$-coordinates for the branch tips of $P_n$ and $(LR)^n$ will be the same, as can be shown by induction on $n$.  Now let $n$ go to infinity to see that the branch tip for the path $\lim_{n \rightarrow \infty} P_n$ will have the same $y$-coordinate as that for $(LR)^\infty$. There are infinitely many such paths. This is why many symmetric binary trees appear to have a ``flat" top as in Figures 2,3, and 4, even when the trees are self-overlapping. As observed, however, by Don West [4], the branch tip for $(LR)^\infty$ or $(RL)^\infty$ may not always have the largest $y$-coordinate. We will show below that for every $\theta$ that is not of the form $\frac{360^\circ}{k}$ for some integer $k$, there will be a critical value $r_T$ such that if $r > r_T$, then there is a branch tip for that symmetric binary tree with a larger $y$-coordinate. These scaling factors will be fairly large, however, and the resulting binary trees will have massively overlapping branches.

Let $\alpha$ be the complex number $\alpha = e^{\theta i}=\cos{\theta} + i \sin{\theta}$. We can represent the trunk of the tree by the complex number $i$ (which we can also think of geometrically as a vector). Recall that multiplying a vector by $r$ will scale the vector by $r$, that multiplying the vector by $\alpha$ corresponds to a counterclockwise rotation by $\theta$, and multiplying by $\alpha^{-1}$ corresponds to a clockwise rotation by $\theta$. We also have $\alpha^n = e^{n\theta i}=\cos{(n\theta)} + i \sin{(n\theta)}$ and $\alpha^{-n} = e^{-n\theta i}=\cos{(n\theta)} - i \sin{(n\theta)}$ for all integers $n$.

To find the branch tip of the path $(LR)^\infty$ in terms of complex numbers, we add the vectors corresponding to the branches. These vectors are obtained by multiplying the branch vectors by \(\alpha\) and \(\alpha^{-1}\) in an alternating fashion while also scaling by $r$. The branch tip for this path is therefore located at the point \((x,y)\) corresponding to the complex number
\begin{align*}
i + ir\alpha + &ir^2 + ir^3\alpha + ir^4 + ir^5\alpha + ir^6 + ir^7\alpha + \dots  \\ \\
&= i \left( 1 + r^2 + r^4 + \dots \right) + i r\alpha \left( 1 + r^2 + r^4 + \dots \right) \\ \\
&= i \frac{1}{1-r^2} + i r \alpha\frac{1}{1-r^2} \\ \\
&= i \frac{1+r \alpha}{1-r^2} = iz.
\end{align*}
For this path we thus have 
$$y = \text{Im}(iz) = \text{Re}(z) = \dfrac{1+r\, \cos{(\theta)}}{1-r^2}.$$

\begin{figure}[ht] 
   \centering
   \includegraphics[height=2.5in]{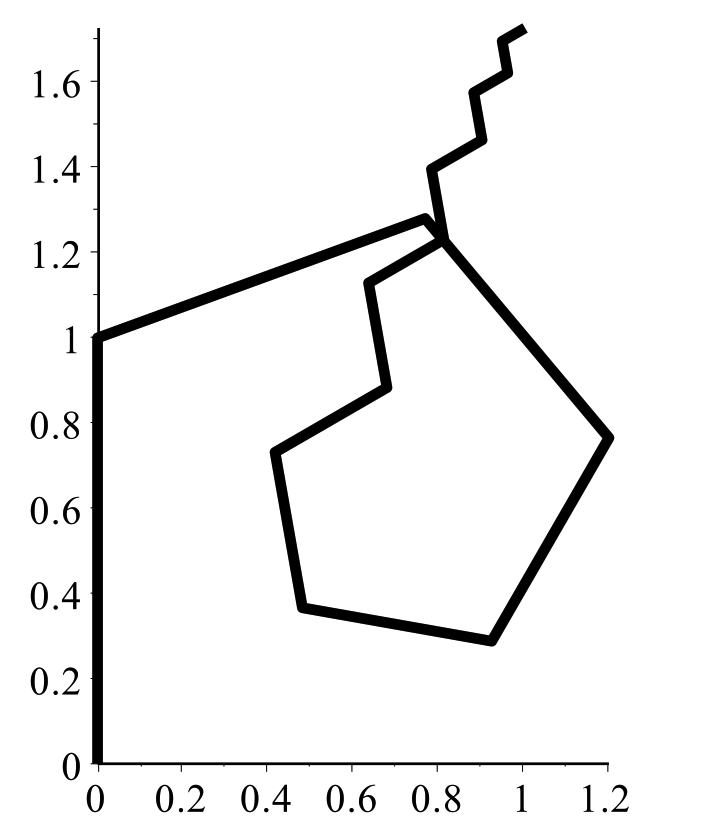} 
   \caption{$R^6(LR)^8$, $\theta = 70^\circ$}
   \label{fig:height example}
\end{figure}

Let $k$ be the smallest integer such that $k\theta \ge 360^\circ$. Consider the branch path that starts $R^k$. This path starts at the top of the trunk and goes clockwise until the last branch regains or passes the upward vertical direction. Now form the rest of the branches by alternating left and right. This yields the branch tip $R^k(LR)^\infty$. Figure~\ref{fig:height example} shows the first 14 branches when $\theta = 70^\circ$ and $k$ = 6.
The branch tip for $R^k(LR)^\infty$ corresponds to the complex number
\begin{align*}
(i + i r \alpha^{-1} &+ i r^2 \alpha^{-2} + \ldots + i r^k \alpha^{-k})  + i r^k \alpha^{-k} \cdot \left( r \alpha + r^2 + r^3 \alpha + r^4 + \dots \right)   \\ \\
&= i \left( \sum_{n=0}^{k}{r^n \alpha^{-n}} + r^k \alpha^{-k} \cdot \left( \dfrac{ra}{1-r^2} + \frac{r^2}{1-r^2} \right) \right)  \\ \\
&= i \left( \sum_{n=0}^{k}{r^n \alpha^{-n}} + \frac{r^{k+1}\alpha^{-(k-1)} + r^{k+2}\alpha^{-k}}{1-r^2} \right) = iz.
\end{align*}
The y-coordinate of this branch tip is \(\text{Im}(iz) = \text{Re}(z)\), or
\[ \sum_{n=0}^{k}{r^n \cos(n \theta)} + \frac{r^{k+1} \cos((k-1)\theta) + r^{k+2} \cos(k \theta)}{1-r^2}.
\]
\subsection{Example}
\begin{figure}[ht] 
   \centering
   \includegraphics[width=5in]{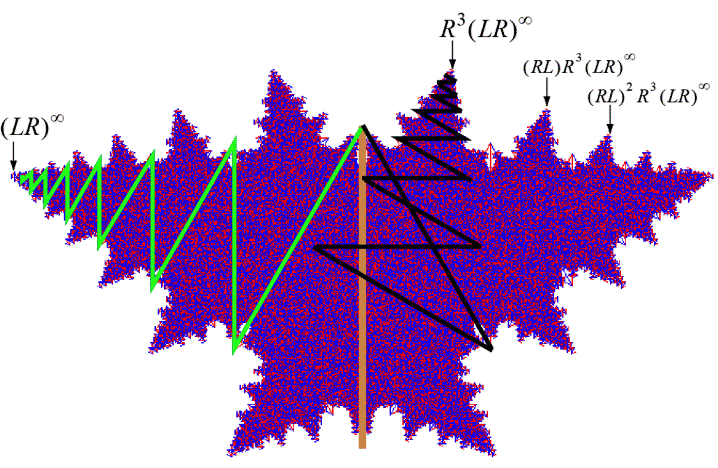} 
   \caption{$\theta = 150^\circ, r= 0.8$}
   \label{fig:150example}
\end{figure}

For a given $\theta$, we are interested in knowing when the branch tip for $R^k(LR)^\infty$ has a greater $y$-coordinate than the branch tip for $(LR)^\infty$. For example, if $\theta = 150^\circ$ and $r = 0.8$, then  Figure~\ref{fig:150example} shows the corresponding symmetric binary tree and the branch tips for these two specific paths (here k = 3). 
The black path for $R^3(LR)^\infty$ has a branch tip with y-coordinate 1.19607 while the green path for $(LR)^\infty$ leads to a branch tip with y-coordinate only at 0.85328. You can get a sense of why this might happen in this example since the green path has a mostly horizontal direction while the black path is more vertical. The value of $r$ is large enough that the black path will grow sufficiently upward to eventually extend above the end of the green path. The branch tip for the $R^3(LR)^\infty$ path corresponds to the top of this tree. Also notice in Figure 6 that the other ``peaks" along the boundary are extensions of $R^3(LR)^\infty$ using initial branches of the form $(RL)^m$. The tree looks "filled-in" because the branches greatly overlap.
\subsection{General Case}
Let $y_k$ be the $y$-coordinate of the $R^k(LR)^\infty$ branch tip and let $y_0$ be the 
$y$-coordinate of the $(LR)^\infty$ branch tip. For a fixed $\theta$, these coordinates depend on the value of $r$. 
Let \(f_\theta (r)=y_k - y_0\) be the difference between the $y$-coordinates of the two branch tips, so
\[
f_\theta (r) = \sum_{n=0}^{k}{r^n \cos(n \theta)} + \frac{r^{k+1} \cos((k-1)\theta) + r^{k+2} \cos(k \theta) -1-r\cos( \theta)}{1-r^2}.
\]
If \(f_\theta (r)\)  is negative then $y_k$ is less than $y_0$, and if \(f_\theta (r)\)  is positive then $y_k$ is greater than $y_0$. We consider two different cases depending on the value of $\theta$.
\subsubsection{Case 1: \(\theta = \dfrac{360^{\circ}}{k}\) where $k$ is an integer}
Note that $k \ge 3$ since $\theta < 180^\circ$, and that $\theta < 1^\circ$ for all $k > 360$. The largest possible angle is $120^\circ$ (for $k = 3$).

\begin{figure}[ht] 
   \centering
   \includegraphics[width=3in]{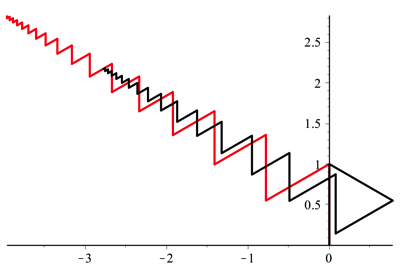} 
   \caption{$\theta = 120^\circ, r= 0.9$}
   \label{fig:120branches}
\end{figure}

For this special case we have $k\theta = 360^\circ$ and therefore the last branch of $R^k$ is vertical, but of length $r^k$. From this point on, the branches of $R^k(LR)^\infty$ alternate in the same pattern that the $(LR)^\infty$ branches did from the beginning, but can never catch up, and thus $y_k$ will be less than $y_0$. Figure~\ref{fig:120branches} shows what happens for $\theta = 120^\circ$ with $r = 0.9$. The red path will have $y_0 = 2.89474$ and the black path will have $y_3 = 2.25526$.

Why does this happen for these special values of $\theta$? Because $k\theta = 360^\circ$, we have \(\cos(k \theta) = \cos(360^{\circ}) = 1\) and \(\cos((k-1)\theta) = \cos(360^{\circ}-\theta) = \cos(\theta)\). Therefore
\[
f_\theta (r) = \sum_{n=0}^{k}{r^n \cos(n \theta)} + \frac{r^{k+1} \cos(\theta) + r^{k+2} -1-r\cos( \theta)}{1-r^2}.
\]

When $r$ is small, the function \(f_\theta (r)\) will behave like the quadratic function $\left( \cos(2\theta)-1 \right)r^2$. Since $0^\circ < \theta < 180^\circ$, the coefficient of $r^2$ in the quadratic will be negative. Therefore the graph of \(f_\theta (r)\) starts at (0,0) and initially begins decreasing in a concave down fashion as $r$ increases. Moreover, with the help of L'Hopital's rule, we can determine that
\begin{align*}
\lim_{r \rightarrow 1^{-}} f_\theta(r) &= \sum_{n=0}^{k}{\cos(n \theta)} + \lim_{r \rightarrow 1^{-}} \frac{(k+1)r^k \cos(\theta) + (k+2)r^{k+1} - \cos(\theta)}{-2r} \\ 
&= 1 + \frac{(k+1)\cos(\theta)+k+2-\cos(\theta)}{-2} = -\frac{k}{2}(1+\cos(\theta)).
\end{align*}
Therefore \(\lim_{r \rightarrow 1^{-}} f_\theta(r) < 0\) suggesting that \(f_\theta(r)\)  remains negative as $r$ goes from 0 to 1. Examining the graphs of \(f_\theta (r)\) for each $k$ from 3 to 360 shows that this is indeed the case for the corresponding $\theta$ values between $1^\circ$ and $120^\circ$. This implies that $y_0$ is always greater than $y_k$ for all $r$ between 0 and 1 for these special values of $\theta$. 
\subsubsection{Case 2: \(\theta \neq \dfrac{360^{\circ}}{n}\) for any integer $n$}
The smallest integer $k$ such that $k\theta > 360^\circ$ is \(k = \left\lceil \frac{360^\circ}{\theta} \right\rceil\), where \(\lceil x \rceil\) is the ceiling function, i.e. the smallest integer greater than or equal to $x$. As in case 1, the value of $k$ can be any integer greater than or equal to 3.

As before, the graph of \(f_\theta (r)\)  starts at (0,0) and initially begins decreasing in a concave down fashion as $r$ increases. This time, however, we have \(\lim_{r \rightarrow 1^{-}} f_\theta (r) = +\infty\). To see why, we just need to verify that the numerator of the fraction in the expression for \(f_\theta (r)\) is positive when $r = 1$ since the denominator is positive and goes to 0 as $r$ approaches 1 from the left. So we need to verify that
\[ N(\theta) = \cos((k-1)\theta) +  \cos(k \theta) -1-\cos( \theta) > 0 \text{ when } \frac{360^\circ}{k} < \theta < \frac{360^\circ}{k-1}
\]
where \(k = \lceil \frac{360^\circ}{\theta} \rceil\). For example, for $120^\circ < \theta < 180^\circ$ with $k = 3$, 
\[N(\theta) = \cos(2\theta) +  \cos(3 \theta) -1-\cos( \theta) = -2 \sin^2(\theta) (2\cos(\theta)+1) > 0.\]
Note that $N(120^\circ) = N(180^\circ) = 0$, and that in general, $N(360^\circ/k) = N(360^\circ/(k-1) = 0$.

\begin{figure}[ht] 
   \centering
   \includegraphics[width=4.5in]{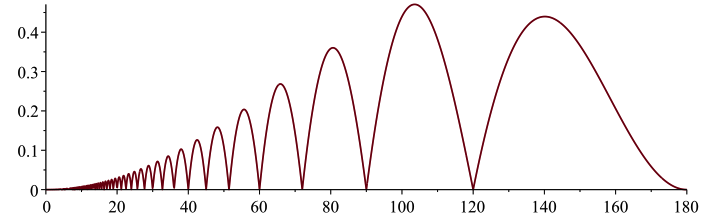} 
   \caption{Graph of $N(\theta)$}
   \label{fig:topnumerator}
\end{figure}

The graph of \(N(\theta)\) is shown in Figure~\ref{fig:topnumerator}, indicating that \(N(\theta)\) is positive for all $\theta$ except for those of the form $\frac{360^\circ}{k}$ (which pile up very quickly as $k$ gets large and $\theta$ approaches 0.)

\begin{figure}[ht]
   \centering
   \begin{subfigure}[b]{0.45\textwidth}
       \includegraphics[width=2.5in]{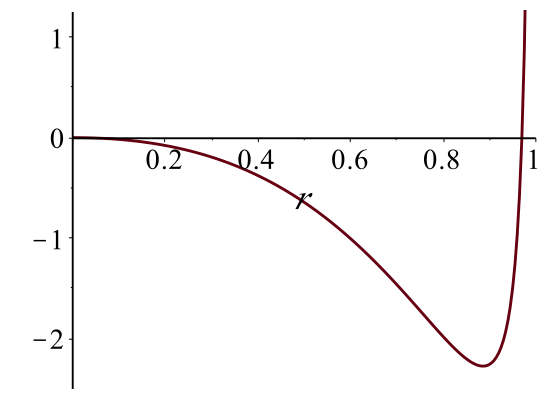}
       \caption{$\theta = 65^\circ, r_T = 0.96984$}
   \end{subfigure}
   \begin{subfigure}[b]{0.45\textwidth}
       \includegraphics[width=2.5in]{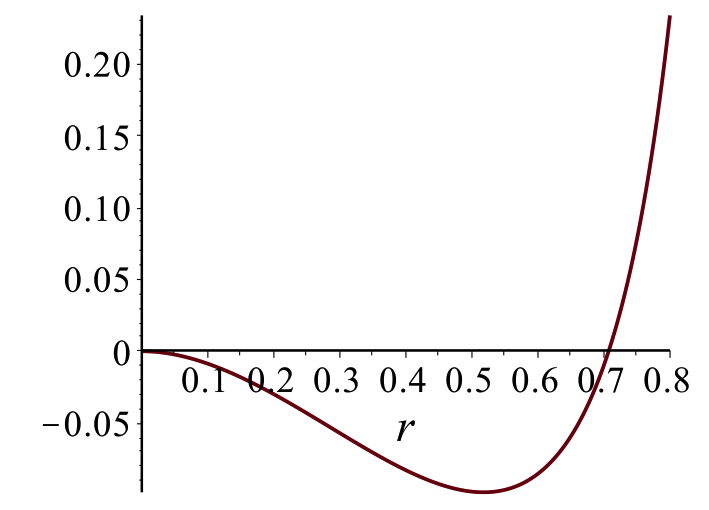}
       \caption{$\theta = 135^\circ, r_T = \frac{1}{\sqrt{2}} = 0.70711$}
   \end{subfigure} 
   \caption{Typical behavior for $f_\theta (r)$ graphs in Case 2}
   \label{fig:topfrgraphs}  
\end{figure}

\begin{figure}[ht] 
   \centering
   \begin{subfigure}[b]{0.45\textwidth}
      \includegraphics[width=\textwidth]{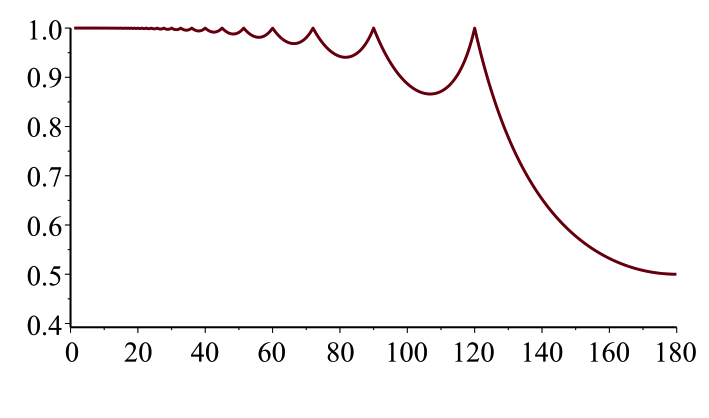} 
      \caption{$0^\circ < \theta < 180^\circ$}
   \end{subfigure}\quad
   \begin{subfigure}[b]{0.45\textwidth}
      \includegraphics[width=\textwidth]{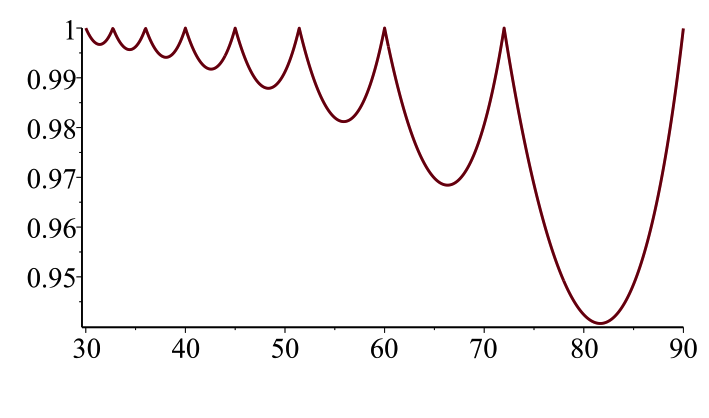} 
      \caption{$30^\circ < \theta < 90^\circ$}
  \end{subfigure}
  \caption{Graphs of $r_T$ as function of $\theta$}
  \label{fig:rTgraphs}
\end{figure}

Because \(f_\theta (r)\) is negative for small $r$ and \(\lim_{r \rightarrow 1^{-}} f_\theta (r) = +\infty\), there must be some some $r_T$ between 0 and 1 where \(f_\theta (r_T) = 0.\) Figure~\ref{fig:topfrgraphs}  shows the graphs of \(f_\theta (r)\) for $\theta = 65^\circ$ and $\theta = 135^\circ$. The shapes of these graphs are typical for all values of $\theta$ in case 2.
The equation \(f_\theta (r) = 0\) can be solved numerically to find the solution $r_T$. Figure~\ref{fig:rTgraphs} shows the graph of $r_T$ as a function of $\theta$ for $0^\circ < \theta < 180^\circ$ and a close up of this graph just on the range $30^\circ < \theta < 90^\circ$.

Several observations can be made. First, the equation \(f_\theta (r) = 0\) has no solution for $r > 0$ when $\theta$ is of the form $\frac{360^\circ}{k}$ for some integer $k$. These are the special angles in case 1. Therefore $r_T$ does not exist at these special values of $\theta$. However,  the graphs show that $r_T$ approaches 1 as $\theta$ approaches one of these special angles.

Second, for $\theta < 90^\circ$, the value of the critical scaling factor $r_T$ is greater than 0.94, so for all practical purposes the top of a reasonable symmetric binary tree for this range of $\theta$ will be determined by the $(LR)^\infty$ path. Any choice of $r > 0.94$ will produce a tree that has a massive overlap of branches.

Third, the situation is a bit different for $\theta > 90^\circ$, and in particular for $\theta > 120^\circ$. On this latter interval, the value of $r_T$ converges to 0.5 as $\theta$ approaches $180^\circ$. We can, in fact, solve \(f_\theta (r) = 0\) explicitly for $r$ for those values of $\theta$ in these two intervals to get
\begin{align*}
90^\circ < \theta < 120^\circ &\Rightarrow r_T =  - \frac{{\cos ( {{{\theta}}} ) + \sqrt { - 3\cos^2 {{( {{{\theta}}} )}} + 1} }}{{4\cos^2 {{( {{{\theta}}} )}} - 1}} \\ 
120^\circ < \theta < 180^\circ &\Rightarrow r_T = -\frac{1}{2\cos(\theta)}
\end{align*}

\begin{figure}[ht] 
   \centering
   \includegraphics[width=3.5in]{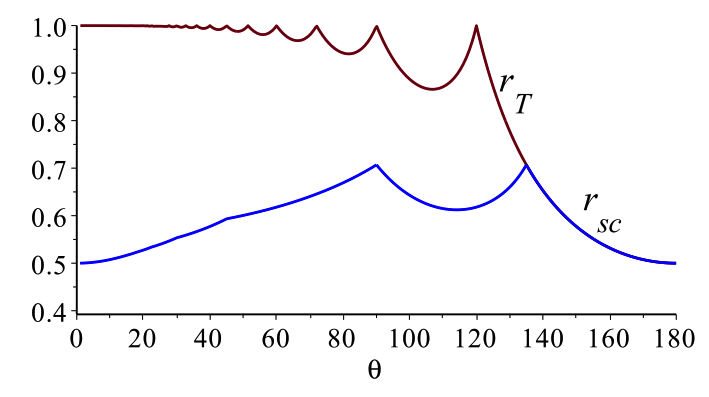} 
   \caption{$r_T \text{ versus } r_{sc}$}
   \label{fig:comparison graphs}
\end{figure}

Finally, it is interesting to note that for $\theta > 135^\circ$ the critical scaling factor $r_T$ is the same value as $r_{sc}$ for which the binary tree of angle $\theta$ is self-contacting. Figure~\ref{fig:comparison graphs} shows a comparison of the two graphs of $r_T$ versus $r_{sc}$. For $r_{sc} < r < r_T$, the binary tree will have overlapping branches, but the top of the branch tips will still correspond to the $(LR)^\infty$ path.

\section{Bottom of the Tree}
\begin{figure}[ht] 
   \centering
   \includegraphics[width=4in]{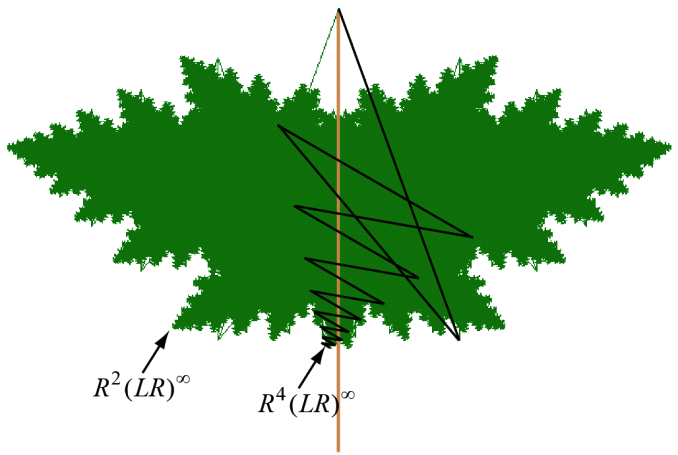} 
   \caption{$\theta = 160^\circ, r = 0.8$}
   \label{fig:SBT160-08}
\end{figure}

As with the top of the tree, we take the ``bottom" of the tree to be the bottom branch tips, that is, the smallest $y$-coordinate of the branch tips. For most trees, the bottom of the tree will correspond to the $y$-coordinate of the branch tip of the path $R^k(LR)^\infty$, where $k$ is now the smallest integer such that $k\theta \ge 180^\circ$. This path bends right until it first achieves or crosses the downward vertical, then alternates left and right. See, for example, Figure~\ref{fig:dimensions}. Just as in the computation in the previous section, the branch tip for the $R^k(LR)^\infty$ path will have $y$-coordinate given by
\[ \sum_{n=0}^{k}{r^n \cos(n \theta)} + \frac{r^{k+1} \cos((k-1)\theta) + r^{k+2} \cos(k \theta)}{1-r^2}.
\]
But as with the top of the tree, however, the branch tip for $R^k(LR)^\infty$, with $k$ the smallest integer satisfying $k\theta \ge 180^\circ$, will not always have the minimum $y$-coordinate. For many angles, if $r$ is sufficiently large, the path $R^m(LR)^\infty$, where $m$ is the smallest integer such that $m\theta \ge 540^\circ$, will have a smaller $y$-coordinate. Since $540^\circ = 180^\circ + 360^\circ$, this path takes another complete set of right turns until it  passes the downward vertical for a second time. For example, this happens for $\theta = 160^\circ$ and $r=0.8$ as illustrated in Figure~\ref{fig:SBT160-08} (here $k = 2$ and $m = 4$). The $R^2(LR)^\infty$ path has a branch tip with $y = 0.27365$ while the black path $R^4(LR)^\infty$ has a branch tip with $y = 0.22498$. 

Let \(f_\theta (r)=y_m - y_k\) be the difference between the $y$-coordinates of the corresponding branch tips for $R^m(LR)^\infty$ and $R^k(LR)^\infty$. Then
$$f_\theta (r) = \sum_{n={k+1}}^{m}{r^n \cos(n \theta)} 
+ \frac{r^{m+1} \cos((m-1)\theta) + r^{m+2} \cos(m \theta) - r^{k+1} \cos((k-1)\theta) - r^{k+2} \cos(k \theta)}{1-r^2}.$$
If \(f_\theta (r)\)  is positive then $y_k$ is less than $y_m$ and if \(f_\theta (r)\)  is negative then $y_k$ is greater than $y_m$. 

Note that \(f_\theta (0) = 0\). For $r$ close to 0, it can be shown that
the graph of \(f_\theta (r)\) will behave either like $-2\sin(\theta)\sin(k\theta)r^{k+1}$ if $k\theta > 180^\circ$, or $(-\cos(2\theta)+1)r^{k+2}$ if $k\theta = 180^\circ$. Suppose $k\theta > 180^\circ$. Since $\theta < 180^\circ$ and $180^\circ < k\theta < 360^\circ$, we have \(\sin(\theta) > 0\) and \(\sin(k\theta) < 0\), and so the coefficient of $r^{k+1}$ in the first case will be positive.  In the second case the coefficient of $r^{k+2}$ will also be positive. Thus in either case the graph of  \(f_\theta (r)\) will initially begin increasing in a concave up fashion as r increases, and thus \(f_\theta (r)\) will be positive for $r$ close to 0. 
\begin{figure}[ht]
   \centering
   \begin{subfigure}[b]{0.3\textwidth}
       \includegraphics[width=\textwidth]{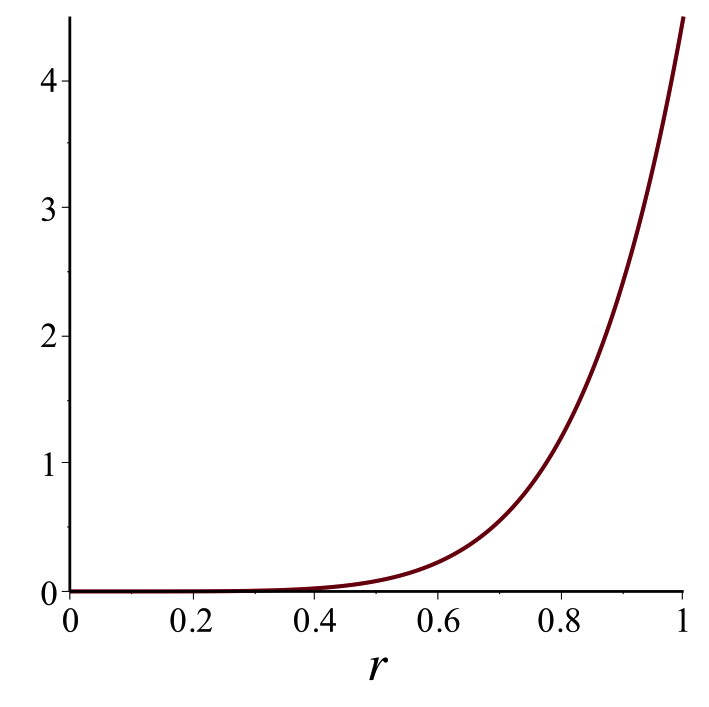}
       \caption{$\displaystyle\lim_{r \rightarrow 1^{-}} f_\theta (r) = A_\theta > 0$}
   \end{subfigure}
   \begin{subfigure}[b]{0.3\textwidth}
       \includegraphics[width=\textwidth]{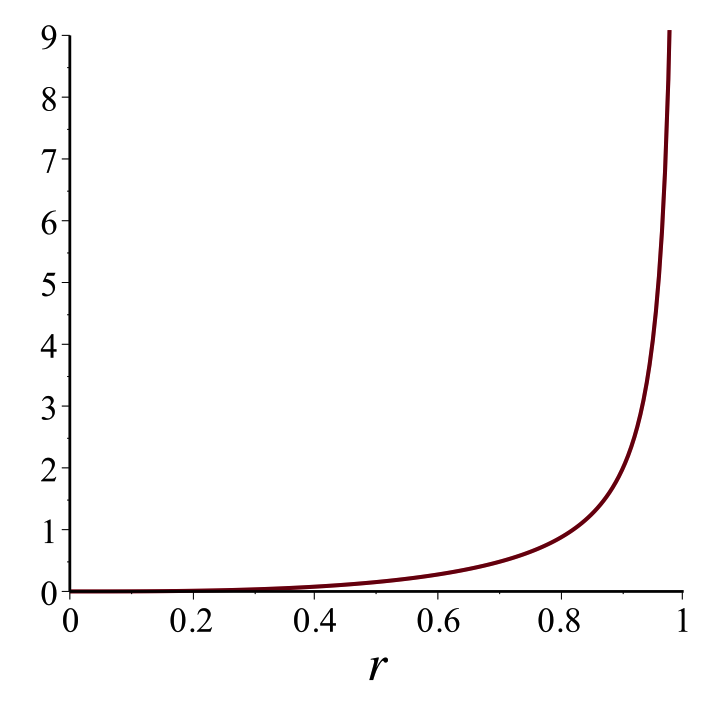}
       \caption{$\displaystyle\lim_{r \rightarrow 1^{-}} f_\theta (r) = +\infty$}
   \end{subfigure}   
   \begin{subfigure}[b]{0.3\textwidth}
       \includegraphics[width=\textwidth]{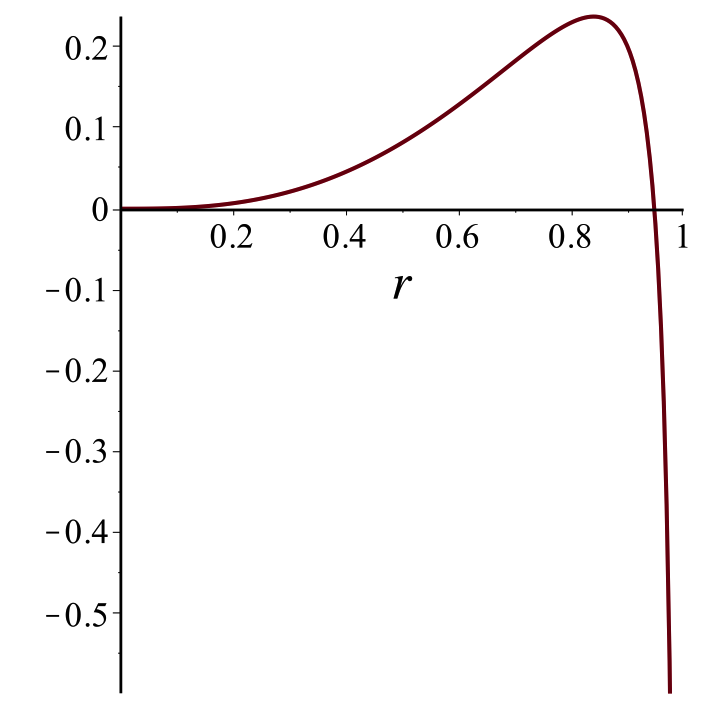}
       \caption{$\displaystyle\lim_{r \rightarrow 1^{-}} f_\theta (r) = -\infty$}
   \end{subfigure}
   \caption{(a) $\theta = 60^\circ$, (b) $\theta = 135^\circ$, (c) $\theta = 146^\circ$ }
   \label{fig:bottom_fr_examples}
\end{figure}

\begin{figure}[ht] 
   \centering
   \includegraphics[width=5in]{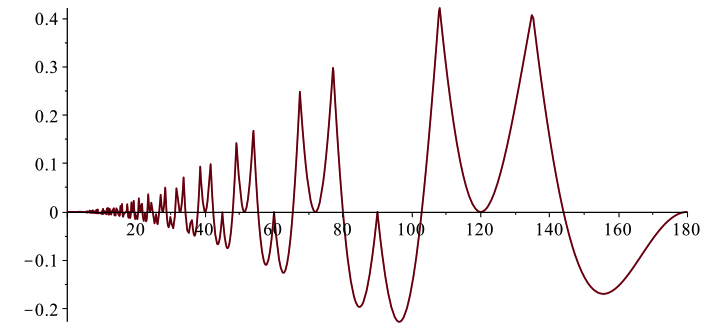} 
   \caption{$N(\theta)$ for $0^\circ < \theta < 180^\circ$}
   \label{fig:bottomNumeratorGraph}
\end{figure}

Figure~\ref{fig:bottom_fr_examples} shows examples of the three things that can happen with the graph of \(f_\theta (r)\) for $0 < r < 1$. Which behavior occurs depends on the sign of the numerator of the fraction in the expression for $f_\theta (r)$ when $r = 1$ since the dominator is positive and goes to 0 as $r$ approaches 1 from the left. So we need to determine the sign of
\[ N(\theta) = \cos((m-1)\theta) +  \cos(m \theta) -\cos((k-1)\theta) -\cos(k \theta)
\]
where \(k = \lceil \frac{180}{\theta}\rceil\) and \(m = \lceil \frac{540}{\theta} \rceil\), in particular whether $N(\theta)$ is 0, positive, or negative (corresponding to graphs (a), (b), and (c), respectively, in Figure~\ref{fig:bottom_fr_examples}). The graph in Figure~\ref{fig:bottomNumeratorGraph} shows that \(N(\theta)\) alternates between intervals where it is positive or zero and intervals where it is negative or zero (with those intervals getting smaller and alternating more rapidly as $\theta$ approaches 0.) 

The places where the graph has a sharp corner are where \(\theta = \frac{540^\circ}{p}\) for an integer $p \ge 4$.
The zeros of $N(\theta)$ occur at \(\theta = \frac{720^\circ}{p}\) for all integers $p \ge 4$. Starting at the right end of the graph and going in decreasing order, the first four zeros are at $180^\circ$, $144^\circ$, $120^\circ$, and $102.85714^\circ$.  

If $\theta$ is a zero of $N(\theta)$, then \(\lim_{r \rightarrow 1^{-}} f_\theta (r)\) exists and is positive. Where $N(\theta)$ is positive we  have \(\lim_{r \rightarrow 1^{-}} f_\theta (r) = +\infty\). For  these values of $\theta$, therefore, \(f_\theta (r) > 0\) for all $0 < r < 1$, and hence $y_k$ will be less than $y_m$.

The intervals where $N(\theta)$ is negative are of the form
\[
\frac{720^\circ}{4n+5} < \theta < \frac{720^\circ}{4n+4} \text{ and } \frac{720^\circ}{4n+4} < \theta < \frac{720^\circ}{4n+3}
\]
for \(n = 0,1,2, \ldots\). (For $n = 0$, only the first interval is used. This gives \(144^\circ < \theta < 180^\circ\).)
On these intervals, we have \(f_\theta (r)\) positive for small $r$ and \(\lim_{r \rightarrow 1^{-}} f_\theta (r) = -\infty\), so there must be some  $r_B$ between 0 and 1 where \(f_\theta (r_B) = 0.\)

\begin{figure}[ht]
   \centering
   \begin{subfigure}[b]{0.45\textwidth}
       \includegraphics[width=\textwidth]{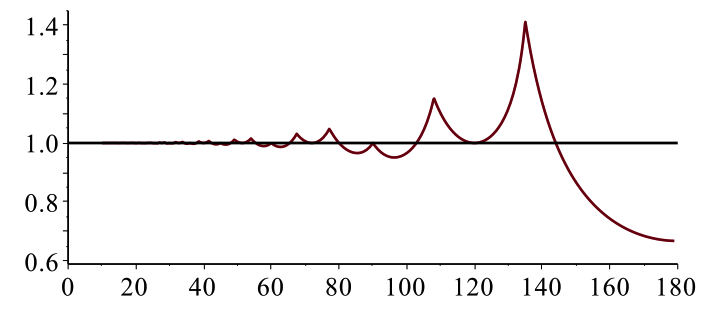}
       \caption{$0^\circ < \theta < 180^\circ$}
   \end{subfigure}
   \begin{subfigure}[b]{0.45\textwidth}
       \includegraphics[width=\textwidth]{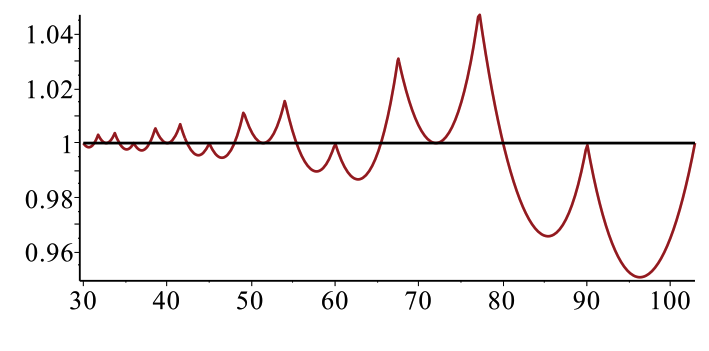}
       \caption{$30^\circ < \theta < 102.85714^\circ$}
   \end{subfigure}
   \caption{Graphs of $r_B$}
   \label{fig:bottomCriticalGraphs}
\end{figure}

The equation \(f_\theta (r) = 0\) can be solved numerically to find the solution $r_B$. This equation will also have a solution on the intervals where $N(\theta)$ is positive, but for those values of $\theta$, $r_B$ will be greater than 1. Figure~\ref{fig:bottomCriticalGraphs} shows the graph of $r_B$ for $0^\circ < \theta < 180^\circ$, and a close up of this graph just on the range $30^\circ < \theta < 102.85714^\circ (= \frac{720^\circ}{7})$.
The intervals where $N(\theta) < 0$ correspond to the  intervals in these graphs where $r_B < 1$. On these intervals, the bottom of the symmetric binary tree will be less than $y_k$ for $r_B < r < 1$. Notice, however, that when $\theta < \frac{720^\circ}{7}$, the value of the critical scaling factor $r_B$ is greater than 0.95 so that for all practical purposes the bottom of a reasonable symmetric binary tree will be determined by the $R^k(LR)^\infty$ path. Any choice of $r > 0.95$ will produce a tree that has a massive overlap of branches. Only on the interval $\theta > 144^\circ$ is there an opportunity for $r_B$ to be small enough that one can generate an interesting tree (such as in Figure~\ref{fig:SBT160-08}). We actually have
\[
r_B = \frac{2\cos(\theta)}{1-4\cos^2 (\theta)}
\]
for $\theta > 144^\circ$ so $r_B$ converges to $\frac{2}{3}$ as $\theta$ approaches $180^\circ$.

\section{Side of the Tree}
Now we want the branch tip that is farthest to the right. Let $k$ be the smallest integer such that $k\theta \ge 90^\circ$. In her thesis, Taylor [4] proved that  the path we seek is $R^k(LR)^\infty$ for self-avoiding and self-contactng trees. This path starts off bending to the right until it first attains or passes a right angle from the initial vertical direction, then alternates left and right branches. The $x$-coordinate for the branch tip for this path will be
\[ \sum_{n=1}^{k}{r^n \sin(n \theta)} + \frac{r^{k+1} \sin((k-1)\theta) + r^{k+2} \sin(k \theta)}{1-r^2}.
\]
By the symmetry of the binary tree, the $x$-coordinate of the branch tip farthest to the left will be the negative of this value. When $\theta \ge 90^\circ$, the largest $x$-coordinate for a self-avoiding or self-contacting tree reduces to
$$x_{\text{max}} = \frac{r\sin(\theta)}{1-r^2}.$$

\begin{figure}[ht]
   \centering
   \begin{subfigure}[b]{0.45\textwidth}
       \includegraphics[width=\textwidth]{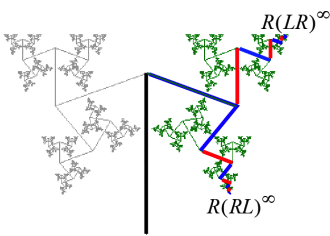}
       \caption{$r =0.60$ (self-avoiding)}
   \end{subfigure}
   \begin{subfigure}[b]{0.45\textwidth}
       \includegraphics[width=\textwidth]{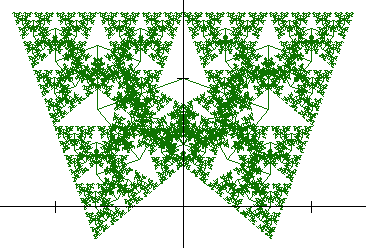}
       \caption{$r =0.71$ (self-overlapping)}
   \end{subfigure}
   \caption{$\theta = 110^\circ$}
   \label{fig:sideExamples}
\end{figure}

When $\theta > 90^\circ$, then $k = 1$ and the path leading to the right side of the tree is $R(LR)^\infty = (RL)^\infty$, which is the same path whose branch tip is usually at the top of the tree. This will happen for all trees with $\theta > 90^\circ$ and $r \le r_{sc}$. Moreover, in this case the bottom tip point will be for the path $R^2(LR)^\infty = R(RL)^\infty$. Thus the side and bottom tip points correspond to the ``corner points" for the scaled subtree whose trunk is the first right branch. Figure~\ref{fig:sideExamples} shows two examples with $\theta = 110^\circ$. The left tree has $r < r_{sc} = 0.61494$, and the right tree has $r > r_{sc}$ but also exhibits this behavior.

\begin{figure}[ht]] 
   \centering
   \includegraphics[width=3in]{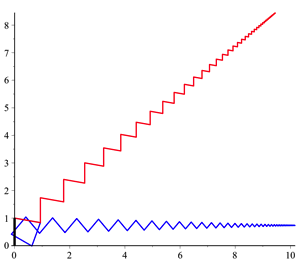} 
   \caption{$\theta = 100^\circ, r = 0.95$ (red: $R(LR)^\infty$, black: $R^5(LR)^\infty$) }
   \label{fig:sidePaths}
\end{figure}

As with the bottom branch tips, however, there are some combinations of angles and scaling factors for which the tree's side will extend past the branch tip described above. The path $R^m(LR)^\infty$, where $m$ is the smallest integer such that $m\theta 
\ge 450^\circ$, may have a larger $x$-coordinate. Since $450^\circ = 90^\circ + 360^\circ$, this path takes another complete set of right turns until it passes the horizontal for the second time. For example, this happens for $\theta = 100^\circ$ and $r=0.95$ as illustrated in Figure~\ref{fig:sidePaths} (here $k = 1$ and $m = 5$). The red path $R(LR)^\infty$ has a branch tip with $x = 9.59556$ while the black path $R^5(LR)^\infty$ has a branch tip with $x = 10.35548$.

The details for which intervals of $\theta$ and values of $r$ the $x$-coordinate of the branch tip for $R^m(LR)^\infty$ is greater than that for the branch tip for $R^k(LR)^\infty$ are similar to those for the bottom of the tree and can be found at Riddle [3]. When $\theta < 72^\circ$, the value of the critical scaling factor for when this happens is greater than 0.99. The other possible intervals are $72^\circ < \theta < 90^\circ$, $90^\circ < \theta < 108^\circ$, and $120^\circ < \theta < 135^\circ$. But even on these three intervals, the critical $r$ value is greater than 0.92 and any such tree would have a massive overlap of branches. So for all practical purposes the right side of a reasonable symmetric binary tree will be determined by the $R^k(LR)^\infty$ path and will have maximum $x$-coordinate given by the formula above.

\section{References}
\begin{enumerate}
\item Mandelbrot, Benoit and Michael Frame. ``The canopy and shortest path of a self-contacting fractal tree," The Mathematical Intelligencer, vol. 21, No. 2 (1999), 18-27.
\item Mandelbrot, Benoit and Michael Frame. ``Geometry of Self-Contacting Binary Fractal Trees," website, \url{http://www.math.union.edu/research/fractaltrees/}.
\item Riddle, Larry. ``The Shape of a Symmetric Binary Tree," website (2014), \url{https://larryriddle.agnesscott.edu/ifs/pythagorean/symbinarytreeShape.htm}.
\item Taylor, Tara. Computational Topology and Fractal Trees, Ph.D. Thesis, Dalhousie University (2005),
\url{http://people.stfx.ca/ttaylor/research/thesis/thesis.pdf}.
\item West, Don. ``Self-Contacting Fractal Trees, Comments on and Mathematical Details for an article by Mandelbrot and Frame," website (1999), \url{http://web.archive.org/web/20120313124012/http://faculty.plattsburgh.edu/don.west/trees/Index.htm}.
\end{enumerate}
\end{document}